\documentclass[12pt]{article}
\usepackage{amssymb, amsfonts, amsmath, amsthm, url, graphicx}

\def\3{\subset }
\def\4{\subseteq }

\def\0{\leqno}

\def\barr{\begin{array}}
\def\earr{\end{array}}

\def\Z{{\rlap{$\kern2pt{\rm Z}$}{\rm Z}\,}}

\title{\bf A characterization of ${\rm PSL}(2,q)$, $q=5,7$}
\author{Marius T\u arn\u auceanu}
\date{February 19, 2016}

\begin{document}

\maketitle

\begin{abstract}
In this short note we prove that the finite non-abelian simple groups ${\rm PSL}(2,q)$,
where $q=5,7$, are determined by their posets of classes of isomorphic subgroups. In particular,
this disproves the conjecture in the end of \cite{5}.
\end{abstract}



\section{Introduction}

In group theory there are many ways to recognize the finite simple groups: by spectrum, by prime graph, by non-commuting graph, by subgroup lattices, ... and so on. Another way to recognize the first two non-abelian simple groups, ${\rm PSL}(2,q)$ with $q=5,7$, is presented in the following. It uses the poset ${\rm Iso}(G)$ of classes of isomorphic subgroups of a group $G$ (see \cite{5}):
$${\rm Iso}(G)=\{[H] \mid H\leq G\}, \mbox{ where }
[H]=\{K\leq G \mid K\cong H\}.$$Recall that ${\rm Iso}(G)$ is partially ordered by
$$[H_1]\leq [H_2] \mbox{ if and only if } K_1\subseteq K_2 \mbox{ for some } K_1\in [H_1] \mbox{ and } K_2\in
[H_2].$$

Obviously, all finite abelian simple groups $G$ have the same poset ${\rm Iso}(G)$  (a chain of length 1) and consequently they cannot be recognized in this way. In the non-abelian case the situation is better, as shows our main result.

\bigskip\noindent{\bf Theorem.} {\it Let $G_0\in\{{\rm PSL}(2,5), {\rm PSL}(2,7)\}$ and $G$ be a finite group
such that ${\rm Iso}(G)\cong{\rm Iso}(G_0)$. Then $G\cong G_0$.}
\bigskip

This leads to a natural question.

\bigskip\noindent{\bf Question.} Let $G_0$ be a finite non-abelian simple group and $G$ be a finite group such that ${\rm Iso}(G)\cong{\rm Iso}(G_0)$. Is it true that $G\cong G_0$?

\section{Proof of the main results}

We start with the following easy but important lemma.

\bigskip\noindent{\bf Lemma.} {\it Let $G_0$, $G$ be two finite groups, $f: {\rm Iso}(G_0)\longrightarrow{\rm Iso}(G)$ be a poset isomorphism and $H_0$, $H$ be two subgroups of $G_0$ and $G$, respectively, such that $f([H_0])=[H]$. Then:
\begin{itemize}
\item[{\rm a)}] ${\rm Iso}(H_0)\cong{\rm Iso}(H)$.
\item[{\rm b)}] If $|H_0|=p_1^{\alpha_1}p_2^{\alpha_2}\cdots p_k^{\alpha_k}$, where $p_i$, $i=1,2,...,k$, are distinct primes, then $|H|=q_1^{\alpha_1}q_2^{\alpha_2}\cdots q_k^{\alpha_k}$ for some distinct primes $q_1,q_2,...,q_k$.
\item[{\rm c)}] If all isomorphic copies of $H_0$ are maximal subgroups of $G_0$, then $H$ is a maximal subgroup of $G$.
\end{itemize}}

\bigskip\noindent{\bf Proof.} a) It is obvious that $f$ induces a poset isomorphism from ${\rm Iso}(H_0)$ to ${\rm Iso}(H)$.

b) This follows from a) and Theorem 3.2 of \cite{5}.

c) Let $K$ be a subgroup of $G$ such that $H\subset K\subset G$. Then $[H]<[K]<[G]$ and therefore $f^{-1}([H])<f^{-1}([K])<f^{-1}([G])$, i.e. $[H_0]<[K_0]<[G_0]$ where $[K_0]=f^{-1}([K])$. It follows that there are $H_0'\cong H_0$ and $K_0'\cong K_0$ such that $H_0'\subset K_0'\subset G_0$, a contradiction.
\hfill\rule{1,5mm}{1,5mm}

\bigskip\noindent{\bf Remark.} The assumption in c) of the above lemma is justified, because a subgroup $M'$ isomorphic to a maximal subgroup $M$ of a group $G$ is not necessarily maximal. For example, let $G$ be a finite non-abelian simple group, $M=\{(x,x)|\,x\in G\}$ and $M'=G\times 1$. Then $M$ is maximal in $G$, $M'\cong M$ $(\cong G)$, but clearly $M'$ is not maximal.
\bigskip

We are now able to prove our main result.

\bigskip\noindent{\bf Proof of the main theorem.} Assume first that $G_0={\rm PSL}(2,5)$. By the above lemma we have $|G|=p^2qr$, where $p$, $q$, $r$ are distinct primes. It is well-known that the maximal subgroups of ${\rm PSL}(2,5)$ are of order $12$ (isomorphic with $A_4$), $10$ (isomorphic with $D_{10}$), and 6 (isomorphic with $S_3$). Moreover, any subgroup of ${\rm PSL}(2,5)$ isomorphic with $A_4$, $D_{10}$ or $S_3$ is also maximal. Therefore, if $f: {\rm Iso}(G_0)\longrightarrow{\rm Iso}(G)$ is a poset isomorphism and $[M_1]=f([A_4])$, $[M_2]=f([D_{10}])$ and $[M_3]=f([S_3])$, then $M_1$, $M_2$ and $M_3$ are maximal subgroups of $G$ of order $p^2q$ (or $p^2r$), $pq$ and $pr$, respectively. Suppose now that $G$ contains a maximal subgroup $M$ which is not isomorphic with $M_1$, $M_2$ or $M_3$. Then $[M]<[M_1]$, $[M]<[M_2]$ or $[M]<[M_3]$ because $[M_1]$, $[M_2]$ and $[M_3]$ are the maximal elements of ${\rm Iso}(G)$. This implies that $|M|$ is a proper divisor of $|M_1|$, $|M_2|$ or $|M_3|$, i.e. $|M|\in\{p,q,r,pq \,({\rm or}\, pr)\}$. Consequently, the orders of maximal subgroups of $G$ are $p^2q$ (or $p^2r$), $pq$, $pr$, and possibly proper divisors of these numbers.

Then $G$ has no Sylow system (it cannot have subgroups of order $qr$) and therefore it is not solvable. Since ${\rm PSL}(2,5)\cong A_5$ is the unique non-solvable group of order $p^2qr$ (see e.g. \cite{1}), it follows that $G\cong G_0$, as desired.
\smallskip

Assume next that $G_0={\rm PSL}(2,7)$. Then $|G|=p^3qr$, where $p$, $q$, $r$ are distinct primes. Since the maximal subgroups of ${\rm PSL}(2,5)$ are of order $24$ (isomorphic with $S_4$) and 21 (isomorphic with the Frobenius group of order 21), a similar argument implies that the orders of maximal subgroups of $G$ are $p^3q$ (or $p^3r$), $qr$, and possibly proper divisors of these numbers.

Again, $G$ has no Sylow system (it cannot have subgroups of order $qr$) and thus it is not solvable. This shows that $G$ has a composition factor, say $G_1/G_2$, which is a non-abelian simple group. Then $|G_1/G_2|$ is even. Moreover, Theorem 1.35 of \cite{3} implies that it is divisible by $4$. This leads to $p=2$, i.e. $|G|=8qr$. Consequently, by \cite{2}, pp. 12-14, we have either $G_1/G_2\cong{\rm PSL}(2,5)$ or $G_1/G_2\cong{\rm PSL}(2,7)$. In the first case we infer that $G$ has order $120$ and so it is isomorphic to one of the following groups: $S_5$, $A_5\times \mathbb{Z}_2$, and ${\rm SL}(2,5)$. This is impossible because all these groups have no maximal subgroup of order $15$. Then $G_1/G_2\cong{\rm PSL}(2,7)$, which shows that $G\cong G_0$. This completes the proof.
\hfill\rule{1,5mm}{1,5mm}

\vspace*{5ex}\small

\hfill
\begin{minipage}[t]{5cm}
Marius T\u arn\u auceanu \\
Faculty of  Mathematics \\
``Al.I. Cuza'' University \\
Ia\c si, Romania \\
e-mail: {\tt tarnauc@uaic.ro}
\end{minipage}

\end{document}